\documentclass[letterpaper]{article}
\usepackage[margin=1in]{geometry}
\usepackage[utf8]{inputenc}
\usepackage{amsmath}
\usepackage{mathtools}
\usepackage{graphicx}
\usepackage{flafter}
\begin{document}

\title{Backstepping feedback control of open channel flow}
\author{Mandy Huo, Sami Malek}
\date{May 16, 2014}
\maketitle

\section{Abstract}

We derive a feedback control law for the control of the downstream flow in a 1-D open channel by manipulating the water flow at an upstream location. We use backstepping for controller design and Lyapunov techniques for stability analysis. Finally, the controller is verified with simulations.

\section{Introduction}

The limitation of global water resources and fluctuations in water needs have sparked interest in the automation of management of water distribution systems. Automation techniques can allow for more efficient operations of open-channel systems by providing, for instance, overflow avoidance and timely supply of desired water flow rate. In this project, we study closed-loop control of 1-D open channel flow. \\

Rabbani derives an open-loop control for tracking a desired flow trajectory in a 1-D open channel \cite{Rabbani}. However, the control in \cite{Rabbani} is derived for systems with the specific initial condition that the initial conditions of the plant and desired trajectory match. The goal of this project is to use feedback control to stabilize the plant to the open-loop equilibrium given any initial condition. Our feedback controller can be combined with Rabbani's feedforward controller to drive any system to any desired flow trajectory. \\

This report is organized as follows. We begin by deriving the equations of the plant which we aim to stabilize. The governing equations for this problem are based on the Hayami model, an appropriate simplification of the Saint-Venant equations. We map this to a convenient PDE-ODE cascade, then proceed with backstepping controller design. Backstepping is useful here for two reasons. First, this technique eliminates sources of instability along the domain when the control is acting only at the boundary \cite{Krsticbndryctrl}, which is the case in our problem. Second, with backstepping we can show that a system is exponentially stable. As we will see in the stability section, this is desirable since our aim is to drive some system to zero. Finally, we present simulations to verify the controller.

\section{Problem formulation}

In this section we introduce the equations of the plant. We begin by describing the physical model, then transforming the system to the form of a PDE-ODE interconnection, a convenient form for backstepping design. 

\subsection{The physical model}
The Saint-Venant equations are commonly used to model unsteady flow in an open channel. In deriving these equations we assume hydrostatic pressure and 1-D flow, that is uniform velocity over the cross-section and horizontal water level across the width of the cross-section. The average channel bed slope as well as the variation in channel width are assumed to be small. Also, boundary friction and turbulence are accounted for through resistance laws such as the Manning-Strickler formula. With these assumptions, the Saint-Venant equations consist of two equations: mass conservation and momentum conservation \cite{Litrico}. The equations are as follows:
\begin{align}
A_{t} + Q_{x} = 0 \label{eq:SVmasscons}\\
Q_{t} + \left(\frac{Q^{2}}{A}\right)_{x} &+ gA(Y_{x} + S_{f} - S_{b}) = 0 \label{eq:SVmomcons}
\end{align}
where $A(x,t)$ is the wetted cross-sectional area (m$^{2}$), $Q(x,t)$ is the discharge (m$^{3}$/s) across $A$, $Y(x,t)$ is the water depth (m), $S_{f}$ is the friction slope (m/m), $S_{b}$ is the bed slope (m/m), and $g$ is the gravitational acceleration (m/s$^{2}$). For the purposes of this report, we make further approximation about the geometry of the problem. We assume a rectangular cross-section and a weir structure positioned downstream at $x=L$ to arrive at the following relations:
\begin{align}
A(x,t) &= B_{0}Y(x,t) \label{eq:assumpA}\\
Z(x,t) &= Y(x,t) + S_{b}(L-x) \label{eq:assumpZ}\\
Q(x,t) &= V(x,t)A(x,t) \label{eq:assumpQ}
\end{align}
where $B_{0}$ is the bed width (m), $Z$ is the absolute water elevation (m), and $V(x,t)$ is the mean water velocity (m/s) across section $A$. We also assume that the weir structure can be modeled by the static relation $Q(L,t) = W(Z(L,t))$, where $W(\cdot)$ is an analytical function. \\

\begin{figure}[h!]
\centering
\includegraphics[scale=1]{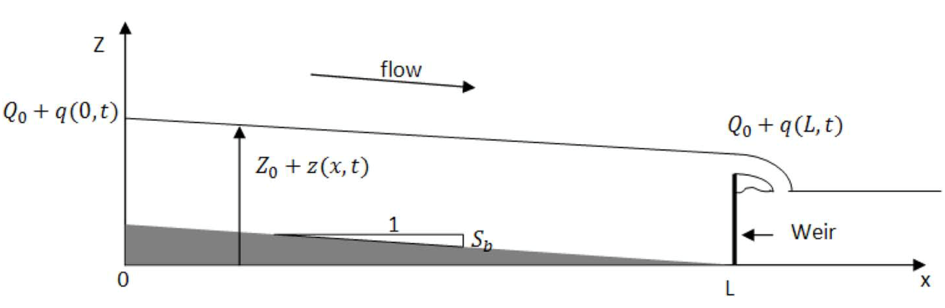}
\caption{1-D cross-section of a water channel with length L.}
\label{fig:waterchann}
\end{figure}

We apply the above relations and linearize the system to arrive at the Hayami model which requires less knowledge of the geometry of the system. The dynamics for $x \in [0,L]$ are modeled by
\begin{align}
D_{0}q_{xx} - C_{0}q_{x} &= q_{t} \label{eq:Hayami1}\\
B_{0}z_{t} + q_{x} &= 0 \label{eq:Hayami2}
\end{align}
where $q(x,t)$ and $z(x,t)$ are the deviations from the nominal flow $Q_{0}$ and the reference elevation $Z_{0}$, respectively, and $C_{0}$ and $D_{0}$ are the nominal celerity and diffusivity, respectively, which depend on $Q_{0}$. A detailed derivation is given in Appendix A of \cite{SVmodelling}. We also assume the relation between the flow and water height at the weir structure can be linearized and thus find the boundary condition
\begin{equation}
q(L,t) = bz(L,t) \label{eq:Hayamibndry}
\end{equation}
where $b$ is the linearization constant (m$^{2}$/s). This is the system treated in \cite{Rabbani}, with the control $u(t) = q(0,t)$ and zero initial conditions. A diagram of the system is shown in Figure \ref{fig:waterchann}.

Our objective is to derive a control which stabilizes the system to the open-loop equilibrium given any (smooth) initial condition. This is equivalent to driving $q(x,t)$ to zero so that $Q(x,t) = Q_{0} + q(x,t) \rightarrow Q_{0}$ over time. We achieve this using backstepping design, but to do so we first map the Hayami model to a system in the form of a PDE-ODE interconnection.

\subsection{The PDE-ODE interconnection}

We study the system approximated by the Hayami model, with the control input $\mu(t)$ actuated at the boundary $x = 0$. We do not use initial conditions in the derivation of the control law, thus the control is general to any (smooth) initial condition. We will make a series of variable changes to make the backstepping design easier. We start by flipping the spatial coordinates so that the input is at $x=1$ and the output is at $x=0$. Applying the variable change $x' = \frac{L-x}{L}$ and renaming the dummy variables, the system \eqref{eq:Hayami1}--\eqref{eq:Hayamibndry} is mapped to
\begin{align} 
z_{t}(x,t)&=\frac{1}{B_{0}}q_{x}(x,t) \label{eq:zq1} \\
q_{t}(x,t)&=D_{0}q_{xx}(x,t)+C_{0}q_{x}(x,t) \label{eq:zq2} \\
q(0,t)&=bz(0,t) \label{eq:zq3} \\
\end{align}
with the control input
\begin{align}
q(1,t)&=\mu(t). \label{eq:zq4}
\end{align}

We set the control law $\mu(t) = e^{\frac{C_{0}}{2D_{0}}x}U(t)$ and define the change of variables $q(x,t) = u(x,t)e^{\frac{C_{0}}{2D_{0}}x}$and $z(x,t) = v(x,t)e^{\frac{C_{0}}{2D_{0}}x}$. Taking the time and spatial derivatives and substituting into \eqref{eq:zq1}--\eqref{eq:zq4}, we have the following system 
\begin{align} 
v_{t}(x,t) &= \frac{1}{B_{0}}u_{x}(x,t) - \frac{C_{0}}{2B_{0}D_{0}}u(x,t) \label{eq:uvsys1} \\
u_t(x,t) &= D_{0}u_{xx}(x,t) - \frac{C_{0}^{2}}{4D_{0}}u(x,t) \label{eq:uvsys2} \\
u(0,t) &= bv(0,t) \label{eq:uvsys3} \\
u(1,t) &= U(t). \label{eq:uvsys4}
\end{align}
\\
Now define $X(t) := B_{0}v(0,t)$. We can then rewrite the $(v,u)$ system in the form of a PDE-ODE interconnection $(X,u)$. We are concerned with the stabilization of the $(X,u)$ system. \\

We differentiate $X(t)$ with respect to time to arrive at the new equation for \eqref{eq:uvsys1} and plug in the variable change to get the new equation for \eqref{eq:uvsys3}. The $(X,u)$ interconnection is

\begin{align} 
\dot{X}(t) &= -\frac{bC_{0}}{2B_{0}D_{0}}X(t) + u_{x}(0,t) \label{eq:Xu1} \\
u_t(x,t) &= D_{0}u_{xx}(x,t) - \frac{C_{0}^{2}}{4D_{0}}u(x,t) \label{eq:Xu2} \\
u(0,t) &= \frac{b}{B_{0}}X(t) \label{eq:Xu3} \\
u(1,t) &= U(t). \label{eq:Xu4}
\end{align}

\section{Controller design}

We design a boundary controller for the PDE-ODE system \eqref{eq:Xu1}--\eqref{eq:Xu4} using backstepping design.

\subsection{Coordinate transformation}

Consider the transformation 
\begin{equation} \label{eq:w}
w(x,t) = u(x,t) - \gamma(x)X(t) - \int^{x}_{0}k(x,y)u(y,t)dy,
\end{equation}
where $\gamma$, $k$, and the control law $U$ are chosen such that \eqref{eq:Xu1}--\eqref{eq:Xu4} is mapped to the following target system
\begin{align} 
\dot{X}(t) &= -\lambda X(t) + w_{x}(0,t) \label{eq:Xw1} \\
w_{t}(x,t) &= D_{0}w_{xx}(x,t) - \frac{C_{0}^{2}}{4D_{0}}w_{x}(x,t) \label{eq:Xw2} \\
w(0,t) &= 0 \label{eq:Xw3} \\
w(1,t) &= 0, \label{eq:Xw4}
\end{align}
where $\lambda$ is some positive constant.
\\ \\
We first derive the kernels $k$ and $\gamma$. The time derivative of the transformation is
\begin{align}
w_{t}(x,t) &= u_{t}(x,t) - \gamma(x)\dot{X}(t) - \int^{x}_{0}k(x,y)u_{t}(y,t)dy \notag \\
&= u_{t}(x,t) - \gamma(x)u_{x}(0,t) + \frac{bC_{0}}{2B_{0}D_{0}}\gamma(x)X(t) - \int^{x}_{0}k(x,y)u_{t}(y,t)dy. 
\end{align}
We want to get rid of the $u_{t}$ term in the integral in order to cancel terms later. Using the PDE on $u$ and integrating by parts twice we get
\begin{align}
\int^{x}_{0}k(x,y)u_{t}(y,t)dy &= D_{0}\int^{x}_{0}k(x,y)u_{yy}(y,t)dy - \frac{C_{0}^{2}}{4D_{0}}\int^{x}_{0}k(x,y)u(y,t)dy \notag \\
&= D_{0} \left[ k(x,x)u_{x}(x,t) - k(x,0)u_{x}(0,t) - k_{y}(x,x)u(x,t) + \frac{b}{B_{0}} k_{y}(x,0)X(t) \right. \notag \\ 
&\quad + \left. \int^{x}_{0} k_{yy}(x,y)u(y,t)dy \right] -\frac{C_{0}^{2}}{4D_{0}}\int^{x}_{0}k(x,y)u(y,t)dy. 
\end{align}
Thus the time derivative is
\begin{align} \label{eq:w_t}
w_{t}(x,t) &= u_{t}(x,t) + \left[ \frac{bC_{0}}{2B_{0}D_{0}}\gamma(x) - \frac{bD_{0}}{B_{0}} k_{y}(x,0)\right]X(t) + \left[D_{0}k(x,0) - \gamma(x)\right] u_{x}(0,t)  \notag \\
&\quad- D_{0} \left[ k(x,x)u_{x}(x,t) - k_{y}(x,x)u(x,t) + \int^{x}_{0} k_{yy}(x,y)u(y,t)dy \right] +\frac{C_{0}^{2}}{4D_{0}}\int^{x}_{0}k(x,y)u(y,t)dy. 
\end{align}

The spatial derivatives are
\begin{align}
w_{x}(x,t) &= u_{x}(x,t) - \gamma'(x)X(t) - k(x,x)u(x,t) - \int^{x}_{0}k_{x}(x,y)u(y,t)dy \label{eq:w_x}\\
w_{xx}(x,t) &= u_{xx}(x,t) - \gamma''(x)X(t) - 2k_{x}(x,x)u(x,t) \notag \\ &- k_y(x,x)u(x,t) - k(x,x)u_{x}(x,t) - \int^{x}_{0}k_{xx}(x,y)u(y,t)dy \label{eq:w_xx}.
\end{align}

We then plug the relevant equations into \eqref{eq:Xu2} and set the coefficients of unwanted terms to zero to arrive at\eqref{eq:Xw2}. We find that $k(x,y)$ must satisfy the following
\begin{align} 
\frac{d}{dx}k(x,x) &= k_{x}(x,x) + k_{y}(x,x) = 0 \label{eq:k1} \\
 k_{xx}(x,x) &= k_{yy}(x,x) \label{eq:k2} \\
k(x,0) &= \frac{1}{D_{0}}\gamma(x) \label{eq:k3}
\end{align}
which is satisfied by
\begin{equation}\label{eq:ksol}
k(x,y) = \frac{1}{D_{0}}\gamma (x-y), \quad x \geq y
\end{equation}
We also find that $\gamma (x)$ must satisfy the following condition: 
\begin{equation} \label{eq:gamma1}
\gamma''(x) = \frac{C_{0}}{2D_{0}^{2}}\left(\frac{C_{0}}{2} - \frac{b}{B_{0}}\right)\gamma(x) + \frac{D_{0}b}{B_{0}}k_{y}(x,0).
\end{equation}
Plugging \eqref{eq:w}  and \eqref{eq:w_x} into \eqref{eq:Xu1} and \eqref{eq:Xu3}, we find the following boundary conditions on $\gamma (x)$: 
\begin{align}
\gamma'(0) &= -\lambda+\frac{C_{0}b}{2D_{0}B_{0}} - \frac{b}{B_{0}}k(0,0) \label{eq:gamma2} \\
\gamma(0) &= \frac{b}{B_{0}}. \label{eq:gamma3} 
\end{align}
Using the kernel equation \eqref{eq:ksol}, we can rewrite this system as a second-order ODE for $\gamma (x)$. The system is
\begin{align}
&\gamma''(x) + \frac{b}{B_{0}}\gamma'(x) - \frac{C_{0}}{2D_{0}^{2}}\left(\frac{C_{0}}{2} - \frac{b}{B_{0}}\right)\gamma(x) = 0 \\
&\gamma'(0) = -\lambda+\frac{C_{0}b}{2D_{0}B_{0}} - \frac{b^{2}}{B_{0}^{2}D_{0}}  \\
&\gamma(0) = \frac{b}{B_{0}}. 
\end{align}
The explicit solution for $\gamma$ will be derived later in the simulation section as it is not needed for proving stability. \\

Using the backstepping transformation to map the boundary conditions of $u$ to $w$ at $x=1$, we find that the control law is
\begin{equation} \label{eq:ctrlU}
U(t) = \gamma(1)X(t) + \int^{1}_{0}k(1,y)u(y,t)dy.
\end{equation}

\subsection{Inverse transformation}
To show stability of the system we will require the inverse transformation. In this section we prove that it exists and find the form of the kernel equations. Let us consider an inverse transformation of the form
\begin{equation} \label{eq:u}
u(x,t) = w(x,t) - \delta (x) X(t) - \int^{x}_{0}l(x,y)w(y,t)dy.
\end{equation}

Proceeding in the same fashion as before we take the time and spatial derivatives, plug into the $(X,w)$ system, and set unwanted coefficients to zero to map to the $(X,u)$ system. The derivatives are
\begin{align}
u_{t}(x,t) &= w_{t}(x,t) + \lambda\delta(x)X(t) + \left[D_{0}l(x,0) - \delta(x)\right] w_{x}(0,t)  \notag \\
&\quad- D_{0} \left[ l(x,x)w_{x}(x,t) - l_{y}(x,x)w(x,t) + \int^{x}_{0} l_{yy}(x,y)w(y,t)dy \right] +\frac{C_{0}^{2}}{4D_{0}}\int^{x}_{0}l(x,y)w(y,t)dy. \label{eq:u_t} \\
u_{x}(x,t) &= w_{x}(x,t) - \delta'(x)X(t) - l(x,x)w(x,t) - \int^{x}_{0}l_{x}(x,y)w(y,t)dy \label{eq:u_x}\\
u_{xx}(x,t) &= w_{xx}(x,t) - \delta''(x)X(t) - 2l_{x}(x,x)u(x,t) \notag \\ & \quad - k_y(x,x)u(x,t) - k(x,x)w_{x}(x,t) - \int^{x}_{0}l_{xx}(x,y)w(y,t)dy \label{eq:u_xx}
\end{align}

We find that the kernel $l(x,y)$ must satisfy
\begin{align} 
\frac{d}{dx}l(x,x) &= l_{x}(x,x) + l_{y}(x,x) = 0 \label{eq:l1} \\
 l_{xx}(x,x) &= l_{yy}(x,x) \label{eq:l2} \\
l(x,0) &= \frac{1}{D_{0}}\delta(x) \label{eq:l3}.
\end{align}
This is of the same form as the conditions on $k(x,y)$, thus a solution to $l(x,y)$ is
\begin{equation}\label{eq:lsol}
l(x,y) = \frac{1}{D_{0}}\delta (x-y), \quad x \geq y.
\end{equation}

The kernel $\delta(x)$ satisfies the following ODE
\begin{align}
\delta''(x) &= \frac{C_{0}^{2}}{4D_{0}^{2}}\delta'(x) - \frac{\lambda}{D_{0}}\delta(x) \label{eq:delta1}\\
\delta'(0) &= -\frac{bC_{0}}{2B_{0}D_{0}} + \lambda \label{eq:delta2}\\
\delta(0) &= -\frac{b}{B_{0}} \label{eq:delta3}.
\end{align}
As this is a second-order linear, homogeneous ODE with constant coefficients, we know that a solution for $\delta(x)$ exists, and thus the inverse transformation exists. Hence, the backstepping transformation is invertible.

\section{Stability analysis}
We now prove exponential stability of the system. Consider the Lyapunov function

\begin{equation} \label{eq:lyap}
V(t) = \frac{1}{2}\|w(x,t)\|^{2} + \frac{c_{1}}{2}\|w_{x}(x,t)\|^{2} + \frac{c_{2}}{2}X(t)^{2},
\end{equation}

where $\|\cdot\|$ denotes the $L_{2}$-norm of a function of $x$. That is,
\begin{equation}
\|w(x,t)\| = \left(\int^{1}_{0}w(x,t)^{2}dx\right)^{1/2}.
\end{equation}

Taking the derivative along the system \eqref{eq:Xw1}--\eqref{eq:Xw4}
\begin{align}
\dot{V}(t) &= \int^{1}_{0}ww_{t}dx + c_{1}\int^{1}_{0}w_{x}w_{xt}dx + c_{2}X\dot{X} \notag \\
&= D_{0}\int^{1}_{0}ww_{xx}dx - \frac{C_{0}^{2}}{4D_{0}}\|w\|^{2} + c_{1}D_{0}\int^{1}_{0}w_{x}w_{xxx}dx -\frac{c_{1}C_{0}^{2}}{4D_{0}}\|w_{x}^{2}\| - c_{2}\lambda X^{2} + c_{2}Xw_{x}(0,t).
\end{align}
Applying Young's inequality twice,
\begin{align}
\dot{V}(t) &\leq - \frac{C_{0}^{2}}{4D_{0}}\|w\|^{2} -\left(\frac{c_{1}C_{0}^{2}}{4D_{0}}+D_{0}\right)\|w_{x}\|^{2}-\frac{\lambda c_{2}}{2}X^{2} - c_{1}D_{0}\|w_{xx}\|^{2} + \frac{c_{2}}{2\lambda}w_{x}(0,t)^{2} \notag \\
&\leq - \frac{C_{0}^{2}}{4D_{0}}\|w\|^{2} -\left(\frac{c_{1}C_{0}^{2}}{4D_{0}}+D_{0} -\frac{c_{2}}{\lambda}\right)\|w_{x}\|^{2}-\frac{\lambda c_{2}}{2}X^{2} +\left(\frac{c_{2}}{\lambda} - c_{1}D_{0}\right)\|w_{xx}\|^{2}.
\end{align}
Setting $c_{1} < \frac{1}{2}$ and $c_{2} = 2D_{0}\lambda c_{1}$, we find
\begin{align}
\dot{V}(t) \leq -\mu V(t),
\end{align} 
where
\begin{equation}
\mu = \min\left\lbrace\frac{C_{0}^{2}}{2D_{0}},2\left[D_{0} + c_{1}\left(\frac{C_{0}^{2}}{4D_{0}} - 2D_{0}\right)\right],2\lambda^{2}D_{0}c_{1}\right\rbrace.
\end{equation}

Thus the Lyapunov function satisfies
\begin{equation} \label{eq:lyapineq}
V(t) \leq e^{-\mu t}V(0)
\end{equation}

Next we will show exponential stability of $u$. Consider the following norms,
\begin{align}
\Gamma_{1}(t) &= \|w(x,t)\|^{2} + \|w_{x}(x,t)\|^{2} + X(t)^{2} \label{eq:Gamma1}\\
\Gamma_{2}(t) &= \|u(x,t)\|^{2} + \|u_{x}(x,t)\|^{2} + X(t)^{2} \label{eq:Gamma2}.
\end{align}

Taking the norm of the backstepping transformation,
\begin{align}
\|w\|^{2} &= \left\| u - \gamma X - \int^{x}_{0}k(x,y)u(y,t)dy \right\|^{2} \notag \\
&\leq 3\left(\|u\|^{2} + \|\gamma\|^{2}X(t)^{2} + \left\|\int^{x}_{0}k(x,y)u(y,t)dy\right\|^{2} \right).
\end{align}

Applying the Cauchy-Schwarz inequality and using the fact that $\gamma(x)$ and $k(x,y)$ are bounded,
\begin{align}
\|w\|^{2} &\leq \alpha_{1}\|u\|^{2} + \alpha_{2}X^{2}, 
\end{align}

where
\begin{align}
\alpha_{1} &= 3\left(1 + \max_{x\in[0,1]}\int^{1}_{0}k(x,y)^{2}dy\right)\\ 
\alpha_{2} &= 3\|\gamma(x)\|^{2}
\end{align}

Similarly, using \eqref{eq:w_x} and solving for $k(x,x)$ with \eqref{eq:ksol} and \eqref{eq:gamma3}, we get
\begin{align}
\|w_{x}\| &= \left\| u_{x} - \gamma' X - k(x,x)u - \int^{x}_{0}k_{x}(x,y)u(y,t)dy \right\|^{2} \notag \\ 
&\leq 4\left(\|u_{x}\|^{2} + \|\gamma'\|^{2}X(t)^{2} + \left(\frac{b}{B_{0}D_{0}}\right)^{2}\|u\|^{2} +  \left\|\int^{x}_{0}k_{x}(x,y)u(y,t)dy\right\|^{2}\right)\notag \\
&\leq \alpha_{3}\|u\|^{2} + \alpha_{4}\|u_{x}\|^{2} + \alpha_{5}X^{2},
\end{align}

where
\begin{align}
\alpha_{3} &= 4\left(\left(\frac{b}{B_{0}D_{0}}\right)^{2} + \max_{x\in[0,1]}\int^{1}_{0}k_{x}(x,y)^{2}dy\right)\\ 
\alpha_{4} &= 4\\
\alpha_{5} &= 4\|\gamma'(x)\|^{2}.
\end{align}

Similarly, using \eqref{eq:u}, \eqref{eq:u_x}, and \eqref{eq:lsol} we can show the following
\begin{align}
\|u\|^{2} &\leq \beta_{1}\|w\|^{2} + \beta_{2}X^{2} \\
\|u_{x}\|^{2} &\leq \beta_{3}\|w\|^{2} + \beta_{4}\|w_{x}\|^{2} + \beta_{5}X^{2},
\end{align}
where
\begin{align}
\beta_{1} &= 3\left(1 + \max_{x\in[0,1]}\int^{1}_{0}l(x,y)^{2}dy\right)\\
\beta_{2} &= 3\|\delta(x)\|^{2}\\
\beta_{3} &= 4\left(\left(\frac{b}{B_{0}D_{0}}\right)^{2} + \max_{x\in[0,1]}\int^{1}_{0}l_{x}(x,y)^{2}dy\right)\\
\beta_{4} &= 4\\
\beta_{5} &= 3\|\delta'(x)\|^{2}.
\end{align}

Thus we have
\begin{align}
\Gamma_{1}(t) &\leq m_{1} \Gamma_{2}(t) \label{eq:G1G2}\\
\Gamma_{2}(t) &\leq m_{2} \Gamma_{1}(t) \label{eq:G2G1},
\end{align}
where
\begin{align}
m_{1} &= \max\{\alpha_{1} + \alpha_{3}, \alpha_{4}, \alpha_{2} + \alpha_{5}\} \\
m_{2} &= \max\{\beta_{1} + \beta_{3}, \beta_{4}, \beta_{2} + \beta_{5}\}.
\end{align}
Since the kernel equations are all continuous along the domain $x \in [0,1]$, they are bounded as long as they are nice at the boundaries. Then the constants $m_{1}$ and $m_{2}$ are finite.  \\

Next we will prove the following inequalities,
\begin{align}
V(t) &\leq m_{3}\Gamma_{1}(t) \label{VleqG1}\\
\Gamma_{1}(t) &\leq m_{4}V(t) \label{G1leqV}.
\end{align}

For the first inequality,
\begin{align}
V(t) &= \frac{1}{2}\|w(x,t)\|^{2} + \frac{c_{1}}{2}\|w_{x}(x,t)\|^{2} + \frac{c_{2}}{2} X(t)^{2} \notag \\
&\leq m_{3}\left(\|w(x,t)\|^{2} + \|w_{x}(x,t)\|^{2} + X(t)^{2}\right) \notag \\
&\leq m_{3}\Gamma_{1},
\end{align}
where
\begin{equation}
m_{3} = \max\{1,c_{1},c_{2}\}.
\end{equation}

For the second inequality,
\begin{align}
\Gamma_{1}(t) &=2\left( \frac{1}{2}\|w(x,t)\|^{2} + \frac{1}{2}\|w_{x}(x,t)\|^{2} + \frac{1}{2} X(t)^{2} \right) \notag \\
&\leq m_{4}\left(\frac{1}{2}\|w(x,t)\|^{2} + \frac{c_{1}}{2}\|w_{x}(x,t)\|^{2} + \frac{c_{2}}{2} X(t)^{2}\right) \notag \\
&\leq m_{4}V(t) 
\end{align}
where 
\begin{equation}
m_{4} = 2\max\{1,c_{1},c_{2}\}.
\end{equation}

Hence using \eqref{eq:lyapineq}, \eqref{VleqG1}, \eqref{G1leqV}, \eqref{eq:G1G2}, and \eqref{eq:G2G1} we have
\begin{align}
\Gamma_{2}(t) \leq R\Gamma_{2}(0)e^{-\mu t},
\end{align}
where
\begin{equation}
R = m_{1}m_{2}m_{3}m_{4}.
\end{equation}

Since $\Gamma _{2}$ is a norm of the $(X,u)$ system, we have shown that the system is indeed exponentially stable. This means $u(x,t) \rightarrow 0$ and $X(t) \rightarrow 0$ as $t \rightarrow \infty$. Since $X(t) = \frac{1}{B_{0}}v(0,t) = \frac{1}{B_{0}}z(0,t)$, we have that $q(x,t) \rightarrow 0$ and $z(0,t) \rightarrow 0$.

\section{Simulations}
We verify our controller by simulations. To do this we first solving the $(X,u)$ system, then map $u$ back to the physical flow $q$ with the variable change $u(x,t) = q(x,t)e^{-\frac{C_{0}}{2D_{0}}x}$. The parameter values used for the following simulations are $b = 1$ m$^{2}$/s, $B_{0} = 7$ m, $C_{0}$ = 20 m/s, $D_{0} = 1800$ m$^{2}$/s, and $L = 1000$ m. These values were used in \cite{Rabbani} to simulate control of the Hayami model. Recall the $(X,u)$ system is 

\begin{align} 
\dot{X}(t) &= -\frac{bC_{0}}{2B_{0}D_{0}}X(t) + u_{x}(0,t) \label{eq:sim1} \\
u_t(x,t) &= D_{0}u_{xx}(x,t) - \frac{C_{0}^{2}}{4D_{0}}u(x,t) \label{eq:sim2} \\
u(0,t) &= \frac{b}{B_{0}}X(t) \label{eq:sim3} \\
u(1,t) &= U(t) = \gamma(1)X(t) + \int^{1}_{0}k(1,y)u(y,t)dy. \label{eq:sim4}
\end{align}

\begin{figure}[h!]
\centering
\includegraphics[scale=0.4]{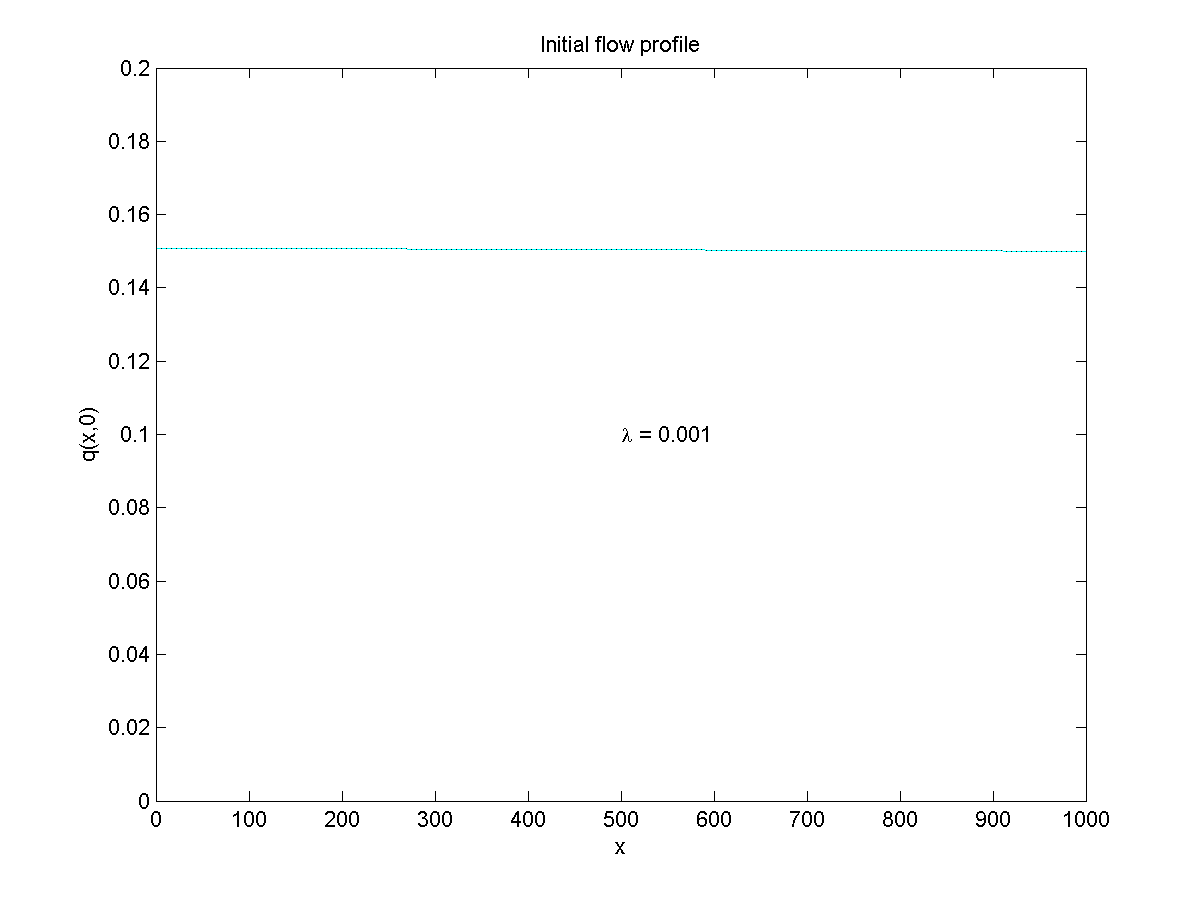}
\caption{Constant initial condition $q(x,t) = 0.15$ m$^{3}$/s. Here $\lambda$ = 0.001}
\label{fig:q_0_1}
\end{figure}

\begin{figure}[h!]
\centering
\includegraphics[scale=0.4]{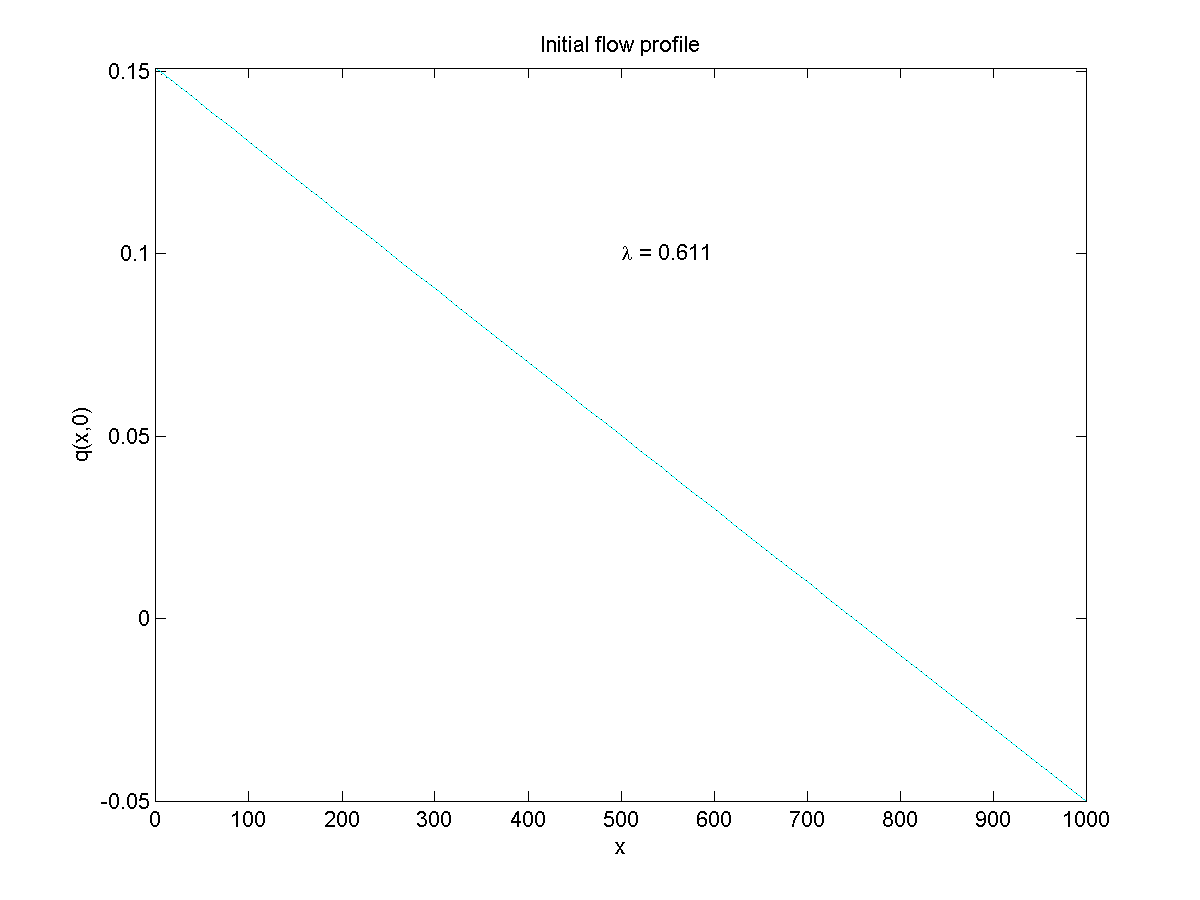}
\caption{Varying initial condition. Here $\lambda$ = 0.611}
\label{fig:q_0_2}
\end{figure}

The control law requires the explicit solutions of the kernel equations for implementation. Earlier we found the equation for the kernel $k(x,y)$ in terms of the kernel $\gamma(x)$. Solving the second-order ODE \eqref{eq:gamma1}--\eqref{eq:gamma3} we find
\begin{align}
\gamma(x) = Ae^{r_{+}x} + Be^{r_{-}x}
\end{align}
where 
\begin{align}
A &= \frac{b}{B_{0}} - B \\
B &= \frac{1}{r_{-} - r_{+}}\left(-\lambda + \frac{C_{0}b}{2D_{0}B_{0}} - \frac{b^{2}}{B_{0}^{2}D_{0}} - \frac{br_{+}}{B_{0}}\right) \\
r_{\pm} &= \frac{-b\pm\sqrt{b^2 + 2\frac{C_{0}B){0}^{2}}{D_{0}^{2}}\left(\frac{C_{0}}{2} - \frac{b}{B_{0}}\right)}}{2B_{0}}.
\end{align}
Note that $\lambda$ was not determined by the backstepping transformation. It remains an arbitrary positive constant which we can tune. Furthermore, the value of $\lambda$ affects $\gamma(x)$ and thus affects the value of the control law $U(t)$. Since the control law is one of the boundary conditions, we choose $\lambda$ so that the control law starts at a value equal to that of the initial condition at the boundary $x=1$, as it should physically. This is done by incrementing $\lambda$ and checking the validity of the control. We consider two arbitrary initial conditions. These are plotted in Figures \ref{fig:q_0_1} and \ref{fig:q_0_2} below along with the initial conditions in which the boundary value at the input is computed from the control law. Note the initial conditions are smooth. The $x$-axis is the spatial coordinate in the physical system, that is $x \in [0, 1000]$.\\

\clearpage

To solve the PDE-ODE system, we approximate the time derivative $u_{t}$ and $\dot{X}$ with a forward difference and $u_{xx}$ with a second-order central difference. We use left-hand rectangle method to discretize the integral when computing $u(1,t)$. This integration method was chosen since it does not require $u(1,t)$ to compute the integral. \\

To deal with the Neumann interconnection, we approximate $u_{x}(0,t)$ with a forward difference so that we can substitute the boundary condition at $x=0$ and arrive at an ODE in terms of only $X(t)$. At each time step $n$, the algorithm is as follows
\begin{enumerate}
\item Solve for $u^{n}_{i}$ except at the boundaries $i = 0$ and $i = end$.
\item Using $u^{n}_{1}$ solve for $X^{n}$. 
\item Compute $u^{n}_{0}$ and $u^{n}_{end}$ using $X^{n}$ and $\gamma$.
\end{enumerate}

The discrete equations are:
\begin{align}
X^{n+1} = \left(1 - \frac{b C_{0} \Delta t}{2B_{0}D_{0}} - \frac{\Delta t}{\Delta x}\frac{b}{B_{0}}\right)X^{n} + \frac{\Delta t}{\Delta x}u_{1}^{n} \\
u^{n+1}_{i} = \frac{\Delta t}{\Delta x}u^{n}_{i+1} + \left(1 - 2\frac{\Delta t}{\Delta x} -\Delta t \frac{C_{0}^{2}}{4D_{0}}\right)u^{n}_{i} + \frac{\Delta t}{\Delta x}u^{n}_{i-1} \\
u_{1}^{n} = \frac{b}{B_{0}}X^{n} \\
u_{x = 1}^{n} = \gamma_{x = 1}X^{n} + \frac{\Delta x}{2D_{0}}\sum^{end - 1}_{i = 0}(\gamma_{end\Delta x - i}u_{i}).
\end{align}

Select time frames of both simulations are show below in Figures \ref{fig:q1} and \ref{fig:q2}. Animations are linked in the Appendix. We see that the control is able to drive the flow deviation $q(x,t)$ to zero, and thus drive the flow to the open-loop equilibrium, even when $q(x,t)$ is non-constant and whether it is initially positive or negative. 
\clearpage
\begin{figure}[h!]
\centering
\includegraphics[scale=0.65]{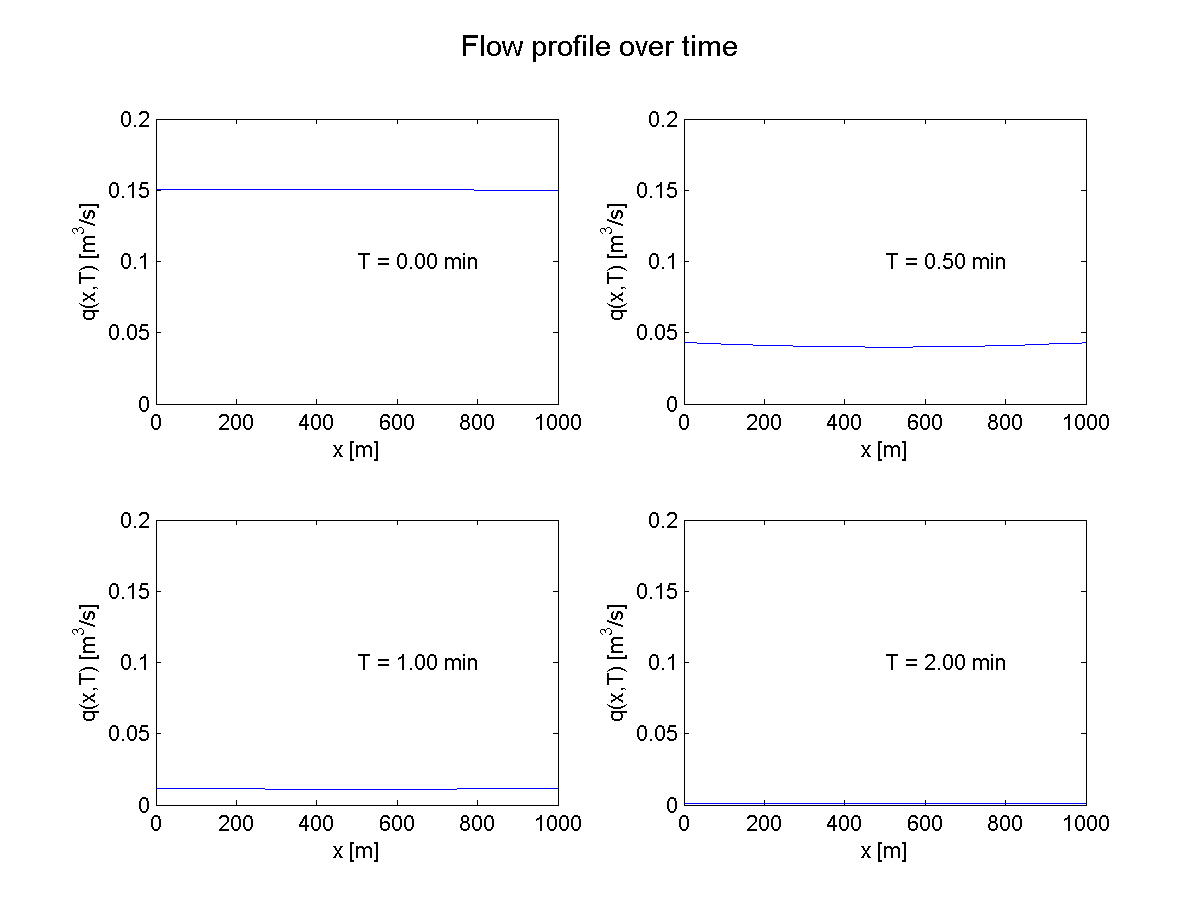}
\caption{Flow profile for constant initial condition $q(x,t) = 0.15$ m$^{3}$/s.}
\label{fig:q1}
\end{figure}
\begin{figure}[h!]
\centering
\includegraphics[scale=0.65]{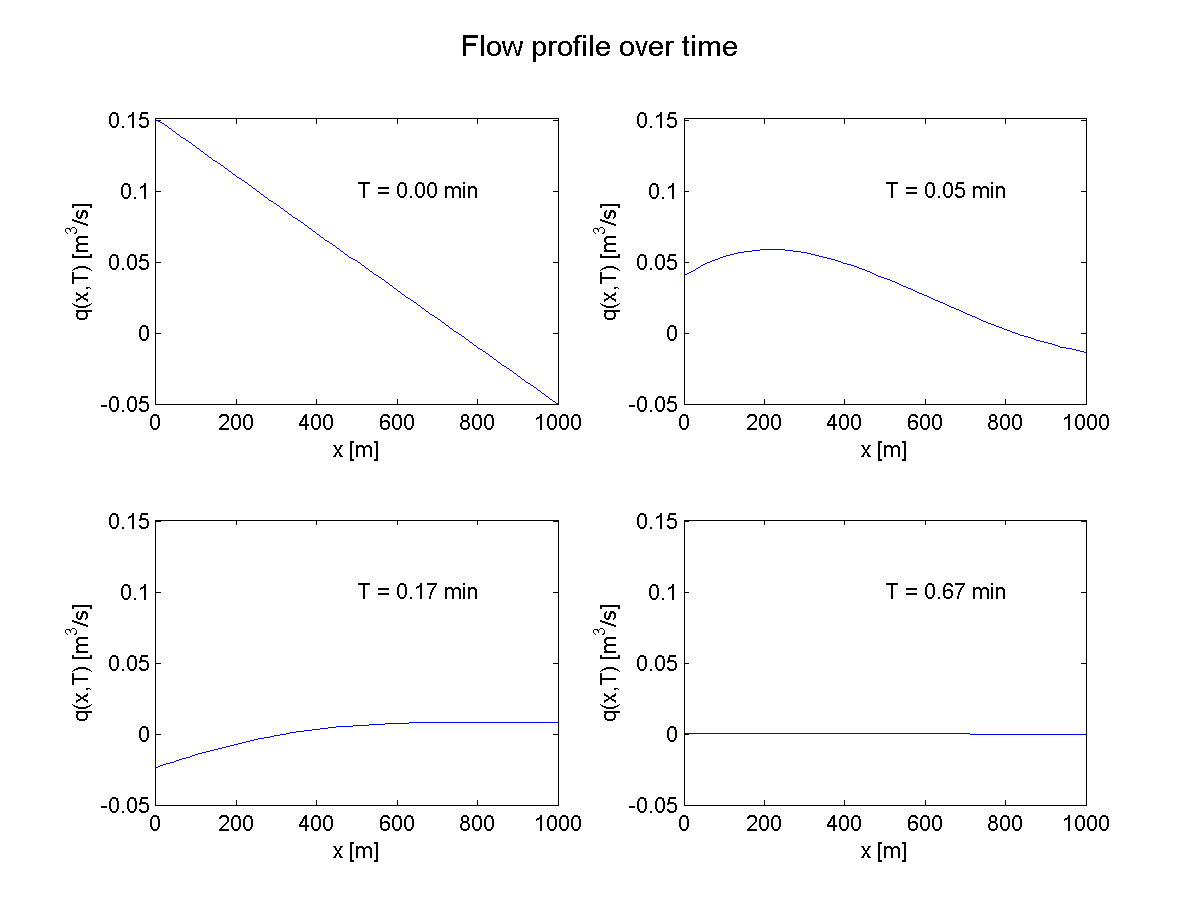}
\caption{Flow profile for non-constant initial condition.}
\label{fig:q2}
\end{figure}

Since $X(t) = \frac{1}{B_{0}}z(0,t)$, where $x = 0$ is the flipped spatial coordinate. In the physical system this gives us the deviation of the water height from the reference $Z_{0}$ at the weir structure. This is plotted below in Figures \ref{fig:q1} and \ref{fig:q2}. We see that the controller also drives $z(x,t)$ to zero. 

\begin{figure}[h!]
\centering
\includegraphics[scale=0.6]{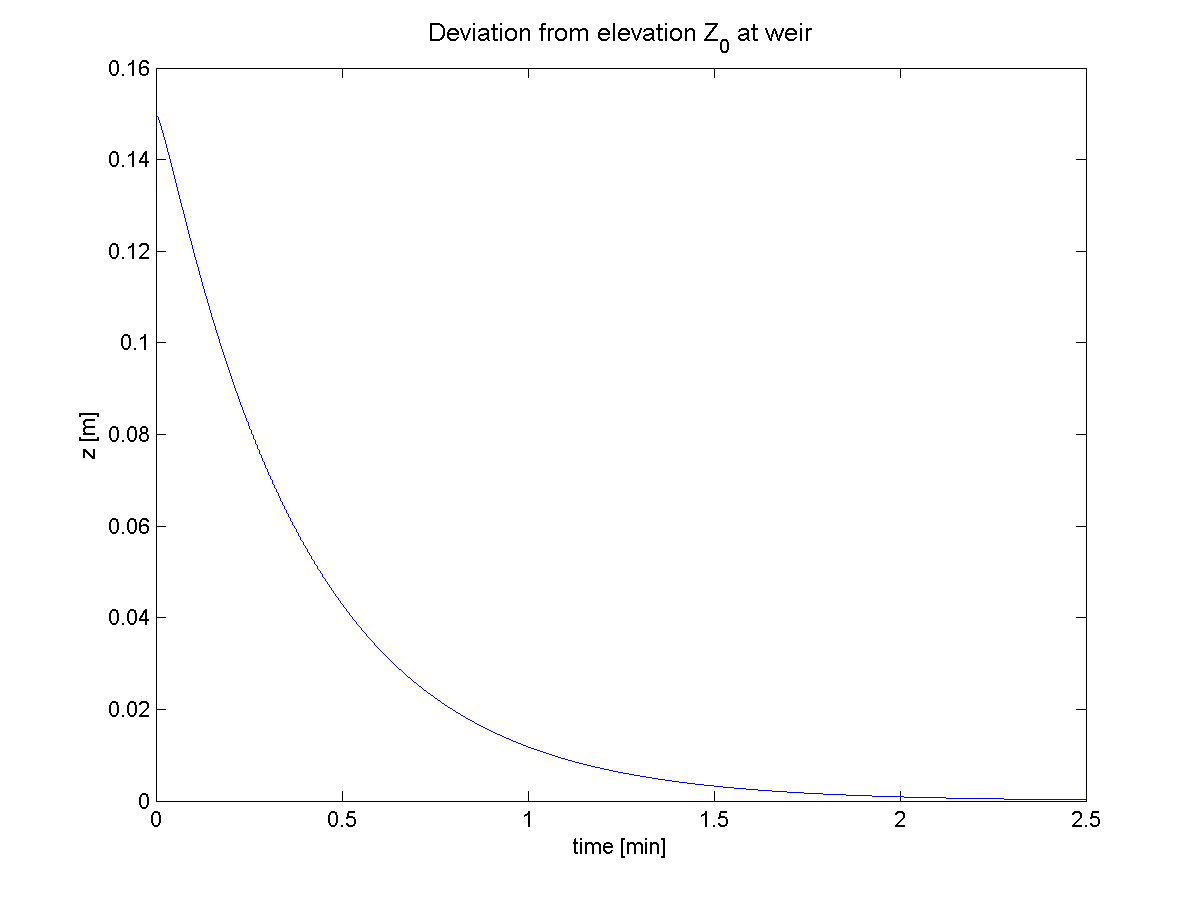}
\caption{Water height deviation at weir over time for constant initial condition.}
\label{fig:z1}
\end{figure}

\begin{figure}[h!]
\centering
\includegraphics[scale=0.6]{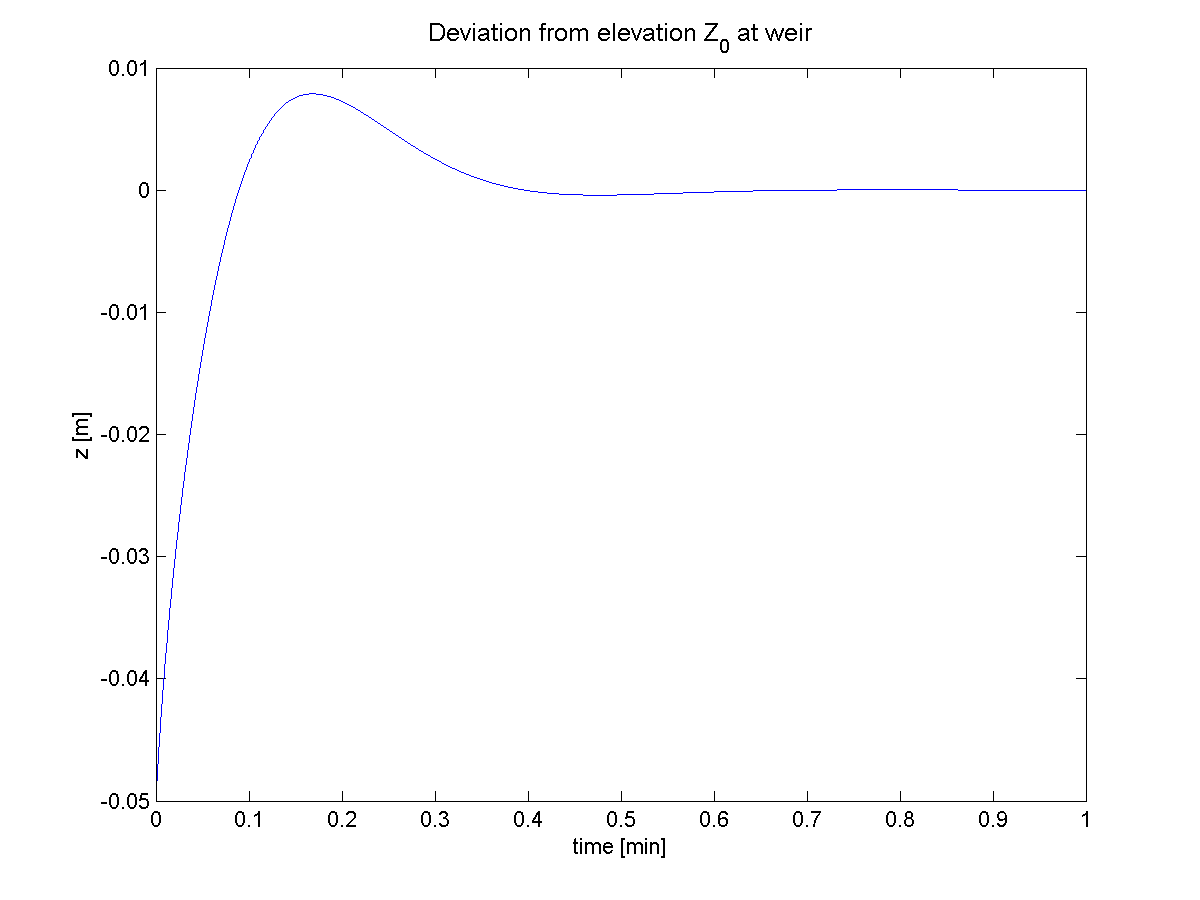}
\caption{Water height deviation at weir over time for non-constant initial condition.}
\label{fig:z2}
\end{figure}

\section{Conclusions}

The controller designed in this report is able to drive any nice initial flow profile to the constant nominal flow $Q_0$. With this feedback control we can build upon the feedforward control designed in \cite{Rabbani}, which works only for specific initial conditions, to derive a control law which allows tracking of a desired flow trajectory under any initial condition. However, it remains to show that the channel will not overflow,that is $z(x,t)$ remains bounded. This can be done by integrating \eqref{eq:uvsys1} in time and using the exponential stability of $u$ to show that $v$ is bounded and thus $z(x,t) = v(x,t)e^{\frac{C_{0}}{2D_{0}}x}$ is also bounded.

\section{Appendix}

The animations for the time evolution of the flow profile can be found at the following links: \\
http://youtu.be/IslywMrgmMY and http://youtu.be/1vFvfBFVp74.

\nocite{*}

\bibliographystyle{ieeetr}
\bibliography{mybibfile}

\end{document}